\documentclass[12pt,reqno]{amsart}
\usepackage{amssymb,amsfonts,latexsym}
\usepackage{graphics,verbatim}
\usepackage{graphicx}
\usepackage[usenames]{color} 
\usepackage{soul} 

\hfuzz=\maxdimen
\hoffset=-1cm \errorcontextlines=0 \numberwithin{equation}{section}

 \theoremstyle{plain}
\newtheorem{theorem}{Theorem}[section]
\newtheorem{lemma}[theorem]{\bf{Lemma}}
\newtheorem{prop}[theorem]{\bf{Proposition}}
\newtheorem{coro}[theorem]{\bf{Corollary}}
\newtheorem{defn}[theorem]{\bf{Definition}}
\newtheorem{remark}[theorem]{\bf{Remark}}

\newenvironment{pf}{{\noindent \bf Proof.\/}}{\hfill$\Box$}

\def \bfi{ {\bf i}}
\def \bfj{ {\bf j}}
\def \bfk{ {\bf k}}

\def \bfw{ {\bf w}}
\def \bfv{ {\bf v}}
\def \bfu{ {\bf u}}

\def \diam{ {\rm diam}}
\date {June, 2012.}

\begin{document}

\title{Graphs induced by iterated function systems}

\author[X.-Y. Wang]{Xiang-Yang Wang}

\address{School of Mathematics and Computational Science\\
Sun Yat-Sen  University\\
Guang-Zhou 510275, P.R. China.}
\email{mcswxy@mail.sysu.edu.cn}


\begin{abstract} For an iterated function system (IFS) of simillitidues, we define two graphs on the representing symbolic space. We show that if the self-similar set $K$ has  positive Lebesgue measure or the IFS satisfies the weak separation condition, then the graphs are hyperbolic, moreover the hyperbolic boundaries are homeomorphic to the self-similar sets.
\end{abstract}


\maketitle

\section{Introduction}

For a contractive iterated function system (IFS) $\{S_j\}_{j=1}^N$ of similitidues on ${\mathbb R}^d$, there is a tree of finite words which represents each point of the associated self-similar set $K$. The iteration defines a random walk on the tree, and the Martin boundary of the random walk is a Cantor set (\cite{[C]}, see also \cite{[K3]}).
On the other hand, Denker and Sato [DS1-3] introduced a random walk on the symbolic space of the Sierpinski gasket, and showed that the Martin boundary is homeomorphic to the gasket. Furthermore, they identified a subclass of harmonic functions from the random walk with Kigami's harmonic functions (\cite{[K1]}, \cite{[K2]}) on the gasket. The case of the pentagasket and other extensions were studied in \cite{[I]} and \cite{[DIK]}. Recently, Lau, Ju and the author \cite{[JLW]} extended this to the class of mono-cyclic post critically finite (p.c.f.) self-similar sets, more generally to self-similar sets with the open set condition (OSC) \cite{[Lau-Wang-2011]}.  This provides a close link of the boundary theory with the recent development of analysis on fractals.

In another direction, Kaimanovich \cite{[Ka]} introduced a hyperbolic structure (``augmented tree") on the symbolic space of the Sierpinski gasket, and showed that the gasket can be identified by the hyperbolic boundary of the graph. The Martin boundary of the simple random walk on the graph can be obtained by a general theory on the random walk on hyperbolic graph (\cite{[A]}, \cite{[Woess]}).

Let $\{S_j\}_{j=1}^N$ be an IFS of similitudes on ${\mathbb R}^d$. Denote by $0< r_i <1, \ i=1,2, \cdots, N$ the contraction ratio of $S_i$. Let $\Sigma^* = \cup_{n=0}^\infty \{1,2,\cdots, N\}^n$ be the finite words space. (We use the notation $o$ to denote the empty word and $\{1,2, \cdots, N\}^0 := \{o\}$). For $\bfi = i_1 \cdots i_n, \bfj = j_1 \cdots j_m \in \Sigma^*$, denote $\bfi \bfj = i_1 \cdots i_n j_1 \cdots j_m$ the concatenation ($o \bfi = \bfi o = \bfi$),  $S_\bfi = S_{i_1} \circ \cdots \circ S_{i_n}$ the composition ($S_o$ is the identity map by convention) and $r_\bfi = r_1 \cdots r_{i_n}$.  Let $r=\min\{r_i:\ i=1,2, \cdots, N\}$, and for each integer $n \ge 1$, let
\[
{\mathcal J}_n = \left\{j_1 j_2 \cdots j_n \in \Sigma^*:\ r_{j_1} \cdots r_{j_n} \le r^n < r_{j_1} \cdots r_{j_{n-1}}\right\}. \quad {\mathcal J} = \cup_{n=0}^\infty {\mathcal J}_n,
\]
where ${\mathcal J}_0=\{o\}$. If $\bfi = i_1 \cdots i_k \in {\mathcal J}_n$, we denote this by $|\bfi| = n$, and say that $\bfi$ is in {\it level} $n$ (Note that $|\bfi|$ is not the length of $\bfi$ in general). We say $\bfi, \bfj \in {\mathcal J}$ are {\it equivalent} and denote by $\bfi \sim \bfj$ if and only if $S_\bfi = S_\bfj$. It is clear that $\bfi \sim \bfj$ implies that $|\bfi|=|\bfj|$. Moreover, $\sim$ defines an equivalence relation on ${\mathcal J}$. We denote by $X$ the quotient space ${\mathcal J}/\sim$, and $[\bfi]$ the equivalence class of $\bfi$. For $x = [\bfi] \in X$, we denote $S_x = S_\bfi$ and $|x| = |\bfi|$. By abusing notation, we write $\bfi \in X$ means that $[\bfi] \in X$.\\

There is a natural graph on $X$: For $x=\{\bfi_1, \cdots, \bfi_n\}, \ y =\{ \bfj_1, \cdots, \bfj_m\} \in X$ (recall an element in $X$ is an equivalence class of some multi-index in $\Sigma^*$), we say that there is an edge between $x$ and $y$ if $\bfj_k = \bfi_\ell \bfk$ for some $1 \le k\le m, 1\le \ell \le n$, $\bfk \in \Sigma^*$ and $|y| = |x| + 1$. We denote by ${\mathcal E}_v$ the above edges set. For $y \in X$, we use the notation $y^{-1}$ to denote any one of $x \in X$ such that $(x,y) \in {\mathcal E}_v$ and $|y| = |x| + 1$. More general, define inductively $y^{-n}= (y^{-(n-1})^{-1}$. It follows that $x, x^{-1}, \cdots, x^{-n}$ is a path from $x$ to $x^{-n}$. If $(x,y) \in {\mathcal E}_v$ with $|y| = |x|+1$, we say that $x$ is an {\it ancestor} of $y$ and $y$ a {\it descendent} of $x$. It is possible that a vertex  in $X$ has more than one ancestor. Also by abusing notation, we write $(\bfi, \bfj) \in {\mathcal E}_v$ to mean that $([\bfi], [\bfj]) \in {\mathcal E}_v$.\\

In order to describe the self-similar set $K$, we need more edges. Let
\[
{\mathcal E}_v^+ = \{ (x,y): \ |y| = |x| + 1,\ S_x(K) \cap S_y(K) \not= \emptyset, \ x \not=y^{-1} \};
\]
and let
\[
{\mathcal E}_h = \{ (x,y): \ |y| = |x| ,\ S_x(K) \cap S_y(K) \not= \emptyset \}.
\]
If $(x,y) \in {\mathcal E}_h$, such that $x^{-1} \not= y^{-1}$ for any $x^{-1}$ and any $y^{-1}$ (recall that $x^{-1}$ may not be unique), then we say $x$ and $y$ are {\it conjugates}. We call an edge in ${\mathcal E}_v \cup {\mathcal E}_v^+$ a {\it vertical edge}, and an edge in ${\mathcal E}_h$ a {\it horizontal edge}. Let
\[
{\mathcal E} = {\mathcal E}_v \cup {\mathcal E}_h, \quad \mbox{and} \quad
{\mathcal E}^\diamond = {\mathcal E}_v \cup {\mathcal E}_v^+.
\]
The graph $(X, {\mathcal E})$ simulates Kaimanovich's ``augmented tree" \cite{[Ka]}.\\

\begin{theorem}\label{theorem-main-1}
The graphs $(X, {\mathcal E})$ and $(X, {\mathcal E}^\diamond)$ are hyperbolic provided that the self-similar set $K$ has positive Lebesgue measure or the IFS satisfies the weak separation condition.
\end{theorem}
The definition of {\it weak separation condition} (WSC) (to see the definition \ref{def-WSC}) was first proposed by Lau and Ngai \cite{[LN]} to study the multifractal structure of an IFS with overlaps, and was studied extensively by many authors (\cite{[Z]}, \cite{[LNR]}, \cite{[FL]}, \cite{[Lau-Wang-2004]} and references there in). The WSC is an important condition in the study of IFS with overlap.

\begin{theorem}\label{theorem-main-2}
With the same assumptions as in Theorem \ref{theorem-main-1}, then the self-similar set $K$ is homeomorphic to the hyperbolic boundaries of $(X, {\mathcal E})$ and $(X, {\mathcal E}^\diamond)$. Furthermore, the H\"{o}lder equivalence holds if we assume additional conditions on the IFS (condition (H) in Section 4).
\end{theorem}

Recall that an IFS satisfies the {\it open set condition} (OSC), if  there exists a bounded nonempty open set $O \subset {\mathbb R}^d$ such that $O \subset \cup_{i=1}^N S_i(O)$ with the union disjoint. It is well known that the OSC implies the WSC, and hence the above theorem \ref{theorem-main-1} and \ref{theorem-main-2} extend the results in \cite{[Lau-Wang-2009]} where the IFS satisfies the OSC.

We organize the paper as following. In Section 2, we recall some basic notations and definitions of a  hyperbolic graph and a hyperbolic boundary.  In Section 3, we study  the properties of the  graphs induced by  an iterated function system, and prove criterions for the graphs $(X, {\mathcal E})$ and $(X, {\mathcal E}^\diamond)$ to be hyperbolic graphs. We prove theorem \ref{theorem-main-1} in Section 4.  In Section 5, we will prove Theorem \ref{theorem-main-2}, and show an example where both the condition (H) and the H\"{o}lder equivalence do not hold. Some open questions are given at the end of the paper.\\

\section{Hyperbolic Graphs and Hyperbolic Boundaries}

Let $G$ be a countably infinite set, and ${\mathcal G} \subset G^2$.  We say that $(G, {\mathcal G})$ (or simply $G$) is a {\it graph} if ${\mathcal G}$ does not have loops and  is symmetric, i.e., $(x,x) \not\in {\mathcal G}$ for all $x\in G$, and $(x,y)\in {\mathcal G}$ implies that $(y,x) \in {\mathcal G}$. We identify $(x,y)$ and $(y,x)$ and call it an {\it edge}. To visualize the graph $(G, {\mathcal G})$, we draw a segment $[x,y]$ if $(x,y) \in {\mathcal G}$. A finite {\it path} $p[x,y]$ from $x$ to $y$ is a sequence $[x_0, x_1, \cdots, x_n]$ with $(x_{i-1}, x_{i}) \in {\mathcal G}$ and $x=x_0, y=x_n$, we use $|p[x,y]| (=n)$ to denote the {\it length} of the path. Throughout the paper, we assume that the graph is {\it connected}, i.e., for any two different vertices $x,y\in G$, there is a path between them. A graph carries an integer-valued metric $d(x,y)$, which is the minimal length of all paths from $x$ to $y$. If a path $p[x,y]$ has the minimal length, we say that the path is a {\it geodesic} segment and denote the path by $\pi[x,y]$. For $x\in G$, we call ${\rm deg}(x) = \{y\in G:\ (x,y)\in {\mathcal G}\}$ the {\it degree} of $x$. We say a graph is {\it local finite} if there exists a constant $c>0$ such that $\max\{{\rm deg}(x):\ x\in G\} \le c$. We fix a reference point $o \in G$ and call it the {\it root}. Denote  $|x|= d(o,x)$, if $|x|=n$, we say $x$ is on the $n$-th {\it level}. If $|x| < |y|$, we say that $x$ is on the {\it upper level } of $y$, or $y$ is on the {\it lower level} of $x$.

Recall that the {\it Gromov product} of two vertices $x,y\in G$ is defined by
\begin{equation}\label{eq-Gromov-def}
|x \wedge y| = \frac{1}{2}(|x| + |y| - d(x,y)).
\end{equation}

\begin{defn} We say a graph $(G, {\mathcal G})$ is $\delta$-hyperbolic (with respect to the root $o$) if there exists a constant $\delta >0$ such that
\begin{equation}\label{def-hyperbolic}
|x \wedge y| \ge \min\{|x \wedge z|, \ |z \wedge y|\} - \delta, \quad \forall x, y, z \in G.
\end{equation}
\end{defn}

As in \cite{[Woess]}, we choose $a>0$ such that $ a' = e^{\delta a}-1
< \sqrt 2 -1$, where $\delta$ is as in (\ref{def-hyperbolic}). Define  for $x,y \in G$,
\begin{equation} \label {eq1.1}
\rho_a(x,y) =
\begin{cases}  \exp(-a|x\wedge y|), & x \not = y,  \\0, \qquad & x=y .
\end{cases}
\end{equation}
 Then
\begin{equation} \label {eq1.2}
\rho_a (x,y) \leq (1+a') \max \{ \rho_a(x,z), \rho_a(y,z)\},
\quad \forall \ x,y,z \in G.
\end{equation}
This means $\rho_a (\cdot , \cdot)$ is an {\it ultra-metric}. It is
not a metric, but is equivalent to the following metric:
$$
\theta_a (x,y) = \inf \{\sum_{i=1}^n \rho_a (x_{i-1},x_i): \ n\geq
1, \ x = x_0, x_1, \dots, x_n =y ,\  x_i\in G\},
$$
in sense that $(1-2 a') \rho_a \le \theta_a \le \rho_a$ (to see \cite[Proposition 22.8]{[Woess]}).
Since $\rho_a$ and $\theta_a$ define the same topology, in our
consideration we will use $\rho_a$ instead of $\theta_a$ for
simplicity. It is known that for any sequence $\{x_n\}_{n=1}^\infty$
such that $\lim_{n\to \infty} |x_n |= \infty$,

\medskip

{\it $\{x_n\}$ is Cauchy in the ultra-metric $\rho_a(x, y)$ if and
only if \ $\lim_{m, n\to \infty} |x_m\wedge x_n| = \infty$.}

\medskip

\begin {defn} Let $\widehat G$ denote the completion of the graph
$G$ under $\rho_a$. We call $\partial G = \widehat G \setminus G$
the {\rm hyperbolic boundary} of $G$.
\end{defn}

\medskip

The hyperbolic boundary  $\partial G$ is a compact set. It  is often
useful to identify $\xi \in  \partial G$ with the {\it geodesic rays} in $G$  that converge to $\xi$. (By a geodesic ray,  we mean an infinite path $ \pi[x_0, x_1, x_2, \cdots]$ such that  $(x_i, x_{i+1}) \in {\mathcal G} \ (i=0, 1, \cdots)$, starting from the root $o$ and with any finite segment of the path being a geodesic). Note that two geodesic rays  $\xi = \pi[x_0, x_1, x_2, \cdots ]$ and $\eta = \pi[y_0, y_1, y_2, \cdots ]$   are equivalent as Cauchy sequences in the ultra-metric $\rho_a$ if and   only if
\begin{equation} \label{eq2.31}
d(x_n, y_n) \leq c\delta
\end{equation}
for all but finitely many $n$, where $c>0$ is independent of the
rays \cite{[Woess]}.\\

Let $\pi[x_0, x_1, x_2, \cdots]$ be a geodesic ray and $y\in G$. For each $n$, there is a geodesic $\pi[y, z_1, \cdots, z_k, x_n]$ connecting $y$ and $x_n$. Note that $p[y, z_1, \cdots, z_k, x_n, x_{n+1}]$ is a path from $y$ to $x_{n+1}$. It follows that $d(y, x_{n+1}) \le d(y, x_n)+1$, this implies that $\{|y \wedge x_n|\}_{n=1}^\infty$ is a non-increasing sequence, and hence $\lim_{n\to\infty} |y \wedge x_n|< \infty$ exists. Similarly, if $\pi[y_0, y_1, y_2, \cdots]$ is another geodesic ray, then  $\lim_{n\to\infty} |x_n \wedge y_n|$ exists and is finite. We extend the Gromov product and ultra-metric to $\partial G$:
\begin{equation}\label{eq-extend-Gromov}
|\xi \wedge \eta| = \inf \{\lim_{n\to\infty} |x_n \wedge y_n |:\ x_n, y_n \in X, \ x_n\longrightarrow \xi,\ y_n \longrightarrow \eta\},
\end{equation}
where the infimum is taken over all geodesic rays $\pi[x_0, x_1, x_2, \cdots]$ and $\pi[y_0, y_1, y_2, \cdots]$ converging to $\xi$ and $\eta$  respectively.
Note that the value of $|x \wedge y|$ has the form $\frac{m}{2}$ (where $m$ is an integer), and hence the infimum is reached by some geodesic rays. Let $\pi[x_0, x_1,x_2, \cdots]$ and $\pi[y_0, y_1, y_2, \cdots]$ be geodesic rays which attain the infimum in ({\ref{eq-extend-Gromov}). Let $\pi[z_0, z_1, z_2, \cdots]$ be another geodesic ray which converges to $\gamma\in \partial G$. By (\ref{def-hyperbolic}), we have
\[
|x_n \wedge y_n| \ge \min\{|x_n \wedge z_n|, \ |z_n \wedge y_n|\} - \delta.
\]
Let $n\to\infty$ and taking the limit, we have
\[
|\xi \wedge \eta| \ge \lim_{n\to\infty} \min\{|x_n \wedge z_n|, \ |z_n \wedge y_n|\} - \delta \ge  \min\{| \xi \wedge \gamma|, \ |\gamma \wedge \eta|\} - \delta.
\]
It follows that
\[
\rho_a(\xi, \eta) \le (1+a') \max\{\rho_a(\xi, \gamma),\ \rho_a(\gamma, \eta)\},
\]
where $a'=e^{a \delta} -1$. We see that the extension of $\rho_a$ is still an ultra-metric on $\partial G$ as in (\ref{eq1.2}).\\

On the other hand, if geodesic rays $\pi[x_1', x_2', \cdots]$ and $\pi[y_1', y_2', \cdots]$ converge to the above $\xi$ and $\eta$ respectively, then $\pi[x_1, x_2, \cdots]$ is equivalent to $\pi[x_1', x_2', \cdots]$, and $\pi[y_1, y_2, \cdots]$ is equivalent to $\pi[y_1', y_2', \cdots]$. By (\ref{eq2.31}), we have
\[
d(x_n,x_n') \le c \delta, \quad d(y_n, y_n') \le c \delta.
\]
It follows that
\[
\begin{aligned}
 \Big{|} |x_n \wedge y_n| - |x_n' \wedge y_n'| \Big{|} & = \frac{1}{2}|d(x_n, y_n) - d(x_n', y_n')| \\
& \le \frac{1}{2} \Big{(} |d(x_n, y_n) - d(x_n, y_n')| + |d(x_n, y_n') - d(x_n', y_n')| \Big{)} \\
& \le \frac{1}{2} \left( d(y_n, y_n') + d(x_n, x_n') \right) \le c \delta.
\end{aligned}
\]
Hence
\begin{equation}\label{eq-ultra-ineq}
\rho_a(\xi,\eta) e^{-a c \delta} \le \lim_{n\to\infty} \rho_a(x_n', y_n')  \le \rho_a(\xi, \eta).
\end{equation}
(The last inequality holds, because $\{x_n\}, \ \{y_n\}$ attain the minimum  in (\ref{eq-extend-Gromov})).
This inequality will be used in section 5.\\

\bigskip

\section{Induced Graphs by  IFS}

Let $\{S_j\}_{j=1}^N$ be an IFS of similitudes on ${\mathbb R}^d$. We sue the notation defined in Section 1 where we defined two graphs $(X, {\mathcal E})$ and $(X, {\mathcal E}^\diamond)$. Let $d(x,y)$ and $d^\diamond(x,y)$ be the graph metrics on $(X, {\mathcal E})$ and $(X, {\mathcal E}^\diamond)$ respectively. We select the empty word $o$ as the root of the graphs, then for any $\bfi \in \Sigma^*$, $|\bfi| = d(o, \bfi)  = d^\diamond (o, \bfi)$ (recall the we abuse the notation $\bfi \in X$ for $[\bfi] \in X$).

If the IFS satisfies the OSC, then the graph $(X, {\mathcal E}_v)$ is a {\it tree} (For any $x\in X$, there exists a unique path from the root to $x$), and this case was studied in \cite{[Lau-Wang-2009]}. If the  OSC does not hold, it is possible that $S_\bfi = S_\bfj$ for deferent $\bfi, \bfj \in \Sigma^*$. Hence there are deferent paths from the root to vertex $[\bfi] = [\bfj] \in X$.

\bigskip

\noindent
{\bf Example 1.}
Let $S_i(x)= \frac{1}{2} (x + i), \ x \in{\mathbb R}, \ i=0,1,2$ be an IFS, the self-similar set is $K=[0,2]$. ${\mathcal J}_2/ \sim = \big\{\{00\}, \{01\}, \{02,10\},\{11\},\{12,20\},\{21\},\{22\}\big\}$.\\

\indent The vertex $\{02, 10\}$ have two ancestors $\{0\}$ and $\{1\}$. $(0, 02), (1,02) \in {\mathcal E}_v$ (abusing the notation).

\begin{center}
\begin{figure}[ht]\label{fig-e}
\centerline{\includegraphics[width=15cm,height=4cm]{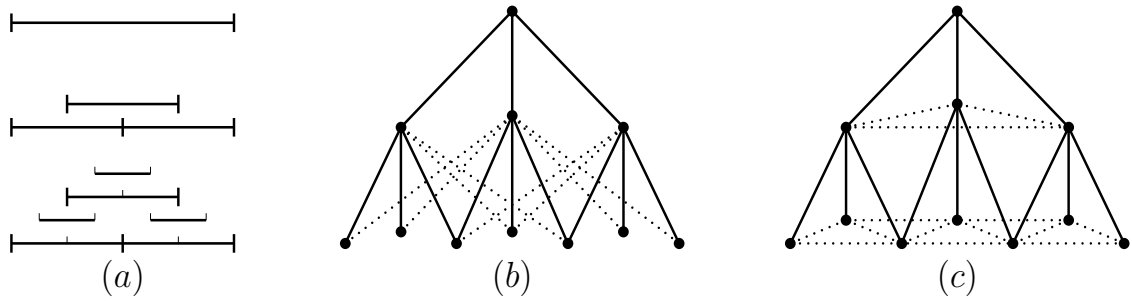}}
\caption{\small{Example 1, $(a)$ the iteration; $(b)$ the graph $(X,{\mathcal E}^\diamond)$; $(c)$ the graph $(X,{\mathcal E})$. The solid lines in $(b), (c)$ are edges in ${\mathcal E}_v$; the doted lines   in $(b), (c)$ are edges in ${\mathcal E}_v^+$, and ${\mathcal E}_h$ respectively.}}
\end{figure}
\end{center}

In the graph $(X, {\mathcal E}^\diamond)$, there are eight edges  connecting $1$: one of them connects the root $o$; three of them connect the descendents ($[10]=\{10, 02\}, [11]=\{11\}, [12]=\{12, 20\}$), and the others belong to ${\mathcal E}_v^+$ (to see Figure 1 $(b)$).

In the graph $(X,{\mathcal E})$, for  $n\ge 2$, each ``boundary vertex" ($0^n$ and $2^n$) has two horizontal neighbors, each ``near boundary vertex" ($0^{n-1}1$ and $2^{n-1}1$) has three horizontal neighbors, and the other vertices have four horizontal neighbors(to see Figure 1 $(c)$).  {\hfill$\Box$}\\

\medskip

For the graph $(X, {\mathcal E})$, a geodesic $\pi[x,y]$ connecting $x$ and $y$ is called {\it canonical} if $\pi[x,y] = \pi[x,u] \cup \pi[u,v] \cup \pi[v,y]$ (one or two parts may vanish) with $\pi[u,v]$ a horizontal path and $\pi[x,u], \ \pi[v,y]$ vertical paths; Moreover for any geodesic path $\pi'[x,y]$, $d(o, \pi[x,y]) \le d(o, \pi'(x,y))$. By the definition of Gromov product (\ref{eq-Gromov-def}), we have
\begin{equation}\label{eq-graph2-Gromov}
|x \wedge y| = h - \frac{\ell}{2},
\end{equation}
where $h$ and $\ell$ are the level and the length of the horizontal segment $\pi[u,v]$ with respectively.\\

Following \cite{[Ka]}, we can use the following moves repeatedly to change the geodesic without increasing  the length: for $u, v \in \pi[x,y]$, \ $|u|=|v|$,
\[
[u, v, v^{-1}] \to [u, u^{-1}, v^{-1}] \quad \mbox{and} \quad [u^{-1}, u, v] \to [u^{-1}, v^{-1}, v].
\]
By using this, we get a canonical geodesic. We should note that for a geodesic segment in ${\mathcal E}$, it cannot contain a sub-segment $[u,v,w]$ with $|u| =|w| = |v| - 1$, since in this case $(u,v), (v,w) \in {\mathcal E}_v$, which implies  that $ S_v(K) \subset S_u(K) \cap S_w(K)$, it follows that $S_u(K) \cap S_w(K) \not= \emptyset$, and hence $(u,w) \in {\mathcal E}_h$ and $d(u,w) = 1$. This contradicts that $[u,v,w]$ is a geodesic segment. \\

\bigskip

An analogous to \cite[Theorem 2.3]{[Lau-Wang-2009]}, we have the following criterion for the graph $(X, {\mathcal E})$ to be hyperbolic. \\

\begin{theorem}\label{theorem-graph-II-hyperbloic}
The graph $(X, {\mathcal E})$ is hyperbolic if and only if there is a constant $L>0$ such that the length of any horizontal geodesic is bounded by $L$.
\end{theorem}

\begin{pf}
The proof of \cite[Theorem 2.3]{[Lau-Wang-2009]}  works here. We give another proof for the necessary part only.

For any horizontal geodesic $\pi[x,y]$ connecting $x, y \in X$, without loss of generality, we assume that the length of $\pi[x,y]$ is an even number, say $2 k$. Let $z$ be the mid-point of $\pi[x,y]$. Then
\[
|x \wedge y| = |x| - k, \ \quad |x \wedge z| = |z \wedge y| = |x| - \frac{k}{2}.
\]
By (\ref{def-hyperbolic}), we have
\[
|x| - k \ge |x| - \frac{k}{2} - \delta, \ \quad {\rm i.e., } \quad k \le 2 \delta.
\]
\end{pf}

\bigskip

Now we study the graph $(X, {\mathcal E}^\diamond)$. Observing that if $p[u,v,w] \ (u\not=w)$ is a path in ${\mathcal E}^\diamond$ such that $|u| = |w| = |v| -1$, ($v$ in the lower level of $u$ and $w$) then $(u, v), \ (v,w) \in {\mathcal E}^\diamond$, and hence $S_u(K) \cap S_v(K) \not= \emptyset$ and $S_v(K) \cap S_w(K) \not= \emptyset$. Let $v' = v^{-2}$ (this $v'$ may not unique), then it is clear that $S_v(K) \subset S_{v'}(K)$. Thus $S_u(K) \cap S_{v'}(K) \not= \emptyset$ and $S_{v'}(K) \cap S_w(K) \not= \emptyset$, it follows that $(u,v'), \ (v', w) \in {\mathcal E}^\diamond$, i.e., $p[u,v',w]$ is also a path. We see that $p[u,v,w,v',u]$ is a closed path with $u, w$ in the same level and $v,\ v'$ in the lower and upper level respectively. The closed path $p[u,v,w,v',u]$ looks like a ``diamond".

\begin{defn} A graph $(G, {\mathcal G})$ is called a {\rm diamond} graph (or simply diamond) if \\
 \indent \hspace{-0.45cm} (i) $(x,y) \not\in {\mathcal G}$ for any $x, y \in G$ with $|x| = |y|$;\\
(ii) For any path $p[u,v,w]$ with $|u| = |w| = |v| -1$, $(u \not= v)$, there exists $v', |v'| = |u|-1$
 such that $p[u,v',w]$ is also a path.\\
\end{defn}

We have shown that

\begin{coro}\label{coro-IFS-diamond}
The graph $(X, {\mathcal E}^\diamond)$ defined in Section 1 is a diamond graph.
\end{coro}

For a diamond graph $(G, {\mathcal G})$, if $[u,v,w]$ is a geodesic segment, then $|v| \not= |u|$ and $|v| \not= |w|$. Hence there are three possible cases: (a) $|u|=|v|+1 = |w|+2$ (or $|w|=|v|+1=|u|+2$); (b) $|u|=|w| = |v|+1$; or (c) $|u|=|w| = |v|-1$. For the last case, we use the move $[u,v,w] \to [u,v',w]$, where $v'$ is as in the above definition. By repeating  this move, we see that for any $x, y \in G$ there is a {\it canonical} geodesic $x= x_0, x_1, \cdots, x_n=y$ such that $|x_i| = |x_{i+1}| +1\ (i < k)$ and $|x_i| = |x_{i+1}| - 1 \ (i \ge k)$ for some $k$, and we say that $x_k$ is on the {\it top level} of the canonical geodesic.
 As a direct consequence of this, we see that $d(x,y)$ is an even number for all $x,y\in G$ with $|x|=|y|$.\\

For a diamond graph, if $x, y \in G$, then there is a canonical geodesic from $x$ to $y$. We assume that $z$ is in the geodesic segment  and is on the top  level. Then it is clear that $|x \wedge y| =|z|$. \\

\begin{theorem}\label{th-hyperbolic-I}
A diamond  graph $(G, {\mathcal G})$ is hyperbolic if and only if there exists some constant $\delta'>0$ such that for any $z\in G$ and any two geodesic pathes $\pi[o, x_1, \cdots, x_n, z]$ and $\pi[o, y_1, \cdots, y_n, z]$ from the root $o$ to $z$, we have $d(x_i, y_i) \le \delta',\ i=1, 2, \cdots, n$.
\end{theorem}

\begin{pf}
{\it Necessity:}
If otherwise, then for any integer $k>0$, there exists $z\in G$ and two geodesic pathes from the root $o$ to $z$: $o\to x \to z$ and $o \to y \to z$, $|x|= |y|$ with $d(x,y) = 2 k$. Let $x= x_1, \cdots, x_k, x^*, x_{k+1}, \cdots, x_{2 k} = y$ be the canonical geodesic joining $x$ and $y$. Then $|x \wedge y| = |x^*| =|x| - k$ and $|x \wedge z| = |z \wedge y| = |x| = |y|$. We see that
\[
|x \wedge y| = \min\{|x \wedge z|, \ |z \wedge y|\} -k.
\]
This contradicts the definition of a hyperbolic graph.\\

\medskip

\noindent
{\it Sufficiency:} We will prove that (\ref{def-hyperbolic}) holds for some  constant $\delta>0$.\\
\indent
For this, we use canonical geodesics connecting  them: $x \to w \to y$, \ $x \to u \to z$ and $z \to v \to y$, where $w, u$ and $v$ are on the top levels, then
\[
|x \wedge y| = |w|, \quad |x \wedge z| =|u|, \quad |z \wedge y|= |v|.
\]

\begin{center}
\begin{figure}[ht]\label{fig-1}
\centerline{\includegraphics[width=5cm,height=3.5cm]{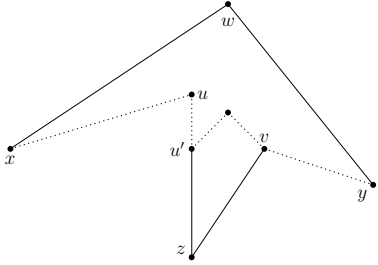}}
\caption{The canonical geodesics}
\end{figure}
\end{center}

Without loss generality, we assume that $|u| \le |v|$. Then (\ref{def-hyperbolic}) is reduced to $|w| \ge |u| - \delta$. Let $u'$ be on the geodesic segment from $u$ to $z$ such that $|u'|=|v|$ (to see Fig 2). The length of the  path from $x$ to $y$: $x \to u \to \ u' \to v \to y$ is $(|x|-|u|) + (|u'|-|u|) + d(u', v) + (|y|-|v|) = |x| + |y| - 2 |u| + d(u', v)$. On the other hand, the canonical geodesic from $x$ to $y$ has length $(|x|-|w|) + (|y|-|w|) = |x|+|y|-2|w|$, and the geodesic has the minimal length. Hence
\[
|x| + |y| - 2 |u| + d(u', v) \ge |x|+|y|-2|w|,
\]
and thus $|w| \ge |u| - \frac{1}{2} d(u', v)$. Now we consider the two geodesics from the root $o$ to $z$: $o \to u' \to z$ and $o \to v \to z$, and note that $|u'|=|v|$. Using the hypothesis, we have $d(u',v) \le \delta'$. It follows that $|w| \ge |u| - \frac{\delta'}{2}$. This completes the proof.
\end{pf}

\bigskip

To end this section, we prove the following lemma which will be used in the next section.

\begin{lemma} \label{lemma-3x}
Let $d$ and $d^\diamond$ be the graph metrics on $(X, {\mathcal E})$ and $(X,{\mathcal E}^\diamond)$ with respectively. Then
\[
d^\diamond(x,y) \le d(x,y) + 1, \quad \forall x, y \in X.
\]
\end{lemma}
\begin{pf}
For any $x,y \in X$, assume that $x, x^{-1}, \cdots, x^{-n}, u_1, \cdots, u_\ell, y^{-m}, \cdots, y^{-1}, y$ is a canonical geodesic in ${\mathcal E}$, where $\{x^{-n}, u_1, \cdots, u_\ell, y^{-m}\}$ is the horizontal part, $\{x, x^{-1}, \cdots, x^{-n}\}$ and $\{y^{-m}, \cdots, y^{-1}, y\}$ are the vertical parts of the canonical geodesic. We consider the two possible cases: (a) $\ell = 2 k +1$ is an odd number; or (b) $\ell = 2 k$ is an even number. In the first case, we replace the horizontal part by $x^{-n}, u_1^{-1}, u_2^{-2}, \cdots, u_k^{-k}, u_{k+1}^{-(k+1)}, u_{k+2}^{-k}, \cdots, u_{2 k +1}^{-1}, y^{-m}$ (this is a path in ${\mathcal E}^\diamond$). Then we get a new path in ${\mathcal E}^\diamond$ with length $\le d(x,y)$.
In the case $\ell = 2 k$, we replace the horizontal part by $x^{-n}, u_1^{-1}, u_2^{-2}, \cdots, u_k^{-k}, u_{k}^{-(k+1)}, u_{k+1}^{-k}, \cdots, u_{2 k}^{-1}, y^{-m}$. We see that $d^\diamond(x,y) \le d(x,y)+1$ in both cases.
\end{pf}

\bigskip

\section{Hyperbolicity of the Graphs}

In this section, we first recall the definition of the {\it weak separation condition} for an IFS and its basic properties. The definition was first proposed by Lau and Ngai \cite{[LN]} to study the multifractal structure of an IFS with overlaps.

\begin{defn}\label{def-WSC}
We say that the IFS $\{S_j\}_{j=1}^N$ satisfies the {\rm weak separation condition} (WSC) if there exists some constant $\gamma>0$ and a compact subset $D \subset {\mathbb R}^d$ with non-empty interior and $\cup_{j=1}^N S_j(D) \subset D$, such that for any $n\ge 1$ and $x\in {\mathbb R}^d$
\[
\# \{ S \in {\mathcal A}_n:\ x \in S(D)\} \le \gamma,
\]
where ${\mathcal A}_n = \{S_\bfi: \ \bfi \in {\mathcal J}_n\}$.
\end{defn}

\begin{lemma}\label{lemma-WSC}
The IFS $\{S_j\}_{j=1}^N$ satisfies the WSC if and only if  for any $b>0$, there exists a constant $\gamma (= \gamma(b))$ such that for any $n$ and $D \subset {\mathbb R}^d$ with $\diam(D) \le b r^n$,
\[
\#\{ x \in X:\ |x|=n,\ S_x(K) \cap D \not= \emptyset\} \le \gamma.
\]
\end{lemma}

This is a consequence of \cite[Proposition 2.1]{[Lau-Wang-2004]}.\\

\begin{theorem}
Assume that the IFS satisfies the weak separation condition. Then the induced graphs $(X, {\mathcal E})$ and $(X, {\mathcal E}^\diamond)$ are local finite.
\end{theorem}

\begin{pf}
For any $x\in X$ with $|x|=n$, let $D= S_x(K)$. Then $\mbox{diam} D \le r^n \mbox{diam} K$. By the above Lemma, we have
\[
\#\{ y\in X: \ |y|= n-1, n \ \mbox{or } n+1,\ S_y(K) \cap D \not= \emptyset \} \le  \gamma(\frac{1}{r}|K|) + \gamma(|K|) + \gamma(r|K|),
\]
where $|K| = \diam(K)$. That graphs are local finite follows from this.
\end{pf}

\bigskip

In the rest of this section, we will prove Theorem \ref{theorem-main-1}. For this, we study the graph $(X, {\mathcal E}^\diamond)$ first.\\

\begin{theorem} \label{theorem-WSC-Hyperbolic}
Suppose the IFS satisfies the WSC, or the self-similar set $K$ has  positive Lebesgue measure. Then the graph $(X, {\mathcal E}^\diamond)$ is hyperbolic.
\end{theorem}

\begin{pf} By Corollary \ref{coro-IFS-diamond}, the graph is diamond. We will make use of Theorem \ref{th-hyperbolic-I} to prove the assertion.

For any $z\in X$ and any two geodesics from the root $o$ to $z$: $o=x_0, x_1, \cdots, x_n = z$; \ $\ o= y_0, y_1, \cdots, y_n = z$. We will prove $d^\diamond(x_k, y_k) \le \delta', \ (k=1,2, \cdots, n)$ for some constant $\delta' >0$ independent of $z$ and the geodesics (where $d^\diamond$ is the graph distance on $(X, {\mathcal E}^\diamond)$). For any fixed $k$,  let
\[
D = \bigcup_{i=k}^n \left( S_{x_i}(K) \cup S_{y_i}(K) \right).
\]
Note that $S_{x_{i-1}}(K) \cap S_{x_i}(K) \not= \emptyset$ and $S_{y_{i-1}}(K) \cap S_{y_i}(K) \not= \emptyset$. It follows that
\[
\begin{aligned}
\mbox{diam} (D) & \le \sum_{i=k}^n \left( \mbox{diam} S_{x_i}(K) + \mbox{diam} S_{y_i}(K) \right)\\
& \le \sum_{i=k}^n 2 r^i |K| < \frac{2 |K|}{1-r} r^k.
\end{aligned}
\]
Let
\[
\{ u_1, \cdots, u_\ell \} = \{ x_i^{-(i-k)}, \ y_i^{-(i-k)}:\ i = k, \ k+1, \cdots, n\}.
\]
Then for each $u_i$, we have $D \cap S_{u_i}(K) \not= \emptyset$, and there exists $\{z_1, \cdots, z_{\ell_0}\} \subset \{u_1, \cdots, u_\ell\}$ such that $x_k = z_1$, $y_k = z_{\ell_0}$ and $S_{z_i}(K) \cap S_{z_{i+1}}(K) \not= \emptyset$ $(i=1, \cdots, \ell_0-1)$. Furthermore, we assume that $\ell_0$ is  minimal. We claim that $\ell_0$ is bounded by some constant.\\

Indeed, if the IFS satisfies the WSC, then by Lemma \ref{lemma-WSC}, we have $\ell_0 \le \gamma(=\gamma(\frac{2|K|}{1-r}))$.\\

Now let us consider the case ${\mathcal L}(K) >0$, where ${\mathcal L}(\cdot)$ is the Lebesgue measure on
${\mathbb R}^d$.  To prove the above $\ell_0$ is bounded by some constant, we let $D' = \cup_{i=1}^{\ell_0} S_{z_i}(K)$. Note that $S_{z_i}(K) \cap D \not= \emptyset$, ${\rm diam}(S_{z_i}(K)) \le r^k {\rm diam}(K)$, $(i=1,2, \cdots, \ell_0)$. Hence
\[
{\rm diam}(D') \le {\rm diam}(D) + 2 r^k {\rm diam}(K) \le \left( \frac{2|K|}{1-r} + 2 |K| \right) r^k := c r^k
\]
By the hypothesis, $\ell_0$ is minimal. We know that each point in $D'$ is covered by at most two $S_{z_i}(K)$. Comparing the Lebesgue measure, we have
\[
r^{(k+1)d} \ell_0 {\mathcal L}(K) \le \sum_{i=1}^{\ell_0} {\mathcal L}(S_{z_i}(K)) \le 2 {\mathcal L}(D')
\le 2 B (c r^k)^d,
\]
where $B$ is the Lebesgue measure of the unite ball in ${\mathbb R}^d$. It follows that $\ell_0 \le \frac{2 B c^d}{r^d{\mathcal L}(K)}$.

This completes the proof of the claim.\\

By the claim, there is a path $x_k = z_1, z_2, \cdots, z_{\ell_0} = y_k$ in $(X, {\mathcal E})$, and hence $d(x_k, y_k) < \ell_0$ (recall that $d$ is the graph metric on $(X, {\mathcal E})$).  By Lemma \ref{lemma-3x}, we have $d^\diamond(x_k, y_k) \le \ell_0$  bounded by some constant.  The assertion follows from this and Theorem \ref{th-hyperbolic-I}.
\end{pf}

\bigskip

In order to prove that the graph $(X, {\mathcal E})$ is hyperbolic, we introduce the following definition.\\

\begin{defn} \label{def-quasi}
Metric space $(X_1, \rho_1)$ is said to be {\rm quasi-isometric} to $(X_2, \rho_2)$
if there exists a map (which is called a quasi-isometry)
$f : X_1 \rightarrow X_2$  and positive constants $\tilde{L}, C$ such that\\
(i) for any $x,y \in X_1$,
\begin{equation} \label{eq-quasi}
\tilde{L}^{-1}\rho_1(x,y) - C < \rho_2(f(x), f(y)) < \tilde{L} \rho_1(x, y) + C;
\end{equation}
(ii) for every $y\in X_2$ there exists $x \in X_1$ such that $\rho_2(y, f(x)) < C$.\\
\end{defn}

Now we can compare the graph metrics $d$ and $d^\diamond$ on $(X, {\mathcal E})$ and $(X, {\mathcal E}^\diamond)$. \\

\begin{theorem}\label{theorem-quasi}
Suppose that the IFS satisfies the WSC or the self-similar set $K$ has positive Lebesgue measure, Then the identity map from the graph $(X, {\mathcal E})$ to $(X, {\mathcal E}^\diamond)$ is a quasi-isometry with the constant $\tilde{L}=1$.
\end{theorem}
\begin{pf}
By Lemma \ref{lemma-3x}, we have $d^\diamond(x, y) \le d(x,y) + 1$.

For the inverse inequality, we assume that $x= x_n, \cdots, x_1, z, y_1, \cdots, y_m = y$ is a canonical geodesic in ${\mathcal E}^\diamond$ with $z$ being on the top level.

Let
\[
D= \left( \cup_{i=1}^n S_{x_i}(K) \right) \bigcup S_z(K) \bigcup \left( \cup_{i=1}^m S_{y_i}(K) \right).
\]
Then
\[
\mbox{diam} D \le \left(\sum_{i=1}^n r^i +1 + \sum_{i=1}^m r^i \right) r^{|z|} |K| < \frac{2 |K|}{1-r} r^{|z|}.
\]

Denote
\[
\{u_1, \cdots, u_\ell\} \subset \{x_n^{-n}, \cdots, x_1^{-1}, z, y_1^{-1}, \cdots, y_m^{-m} \}.
\]
Then $D \cap S_{u_i}(K) \not= \emptyset \ (i=1,\cdots, \ell)$, and exist $\{z_1, \cdots, z_{\ell_0} \} \subset \{u_1, \cdots, u_\ell\}$ such that $z_1 = x_n^{-n}$, $z_{\ell_0} = y_m^{-m}$, and $(z_i, z_{i+1}) \in {\mathcal E}^h$, $(i=1, \cdots, \ell_0 - 1)$. Furthermore, we assume that $\ell_0$ is minimal. Then a similar argument as in the proof of Theorem \ref{theorem-WSC-Hyperbolic} shows that the above $\ell_0$ is bounded by some constant $C>0$. We see that $\pi[x_n, x_n^{-1}, \cdots, x_n^{-n}] \cup p[z_1, \cdots, z_{\ell_0}] \cup \pi[y_m^{-m}, \cdots, y_m^{-1}, y_m]$ is a path from $x= x_n$ to $y = y_m$ in $(X, {\mathcal E})$. Hence
\[
d(x,y) \le n + m + C = d^\diamond(x,y) + C.
\]
This  completes the proof.
\end{pf}

\bigskip

Denote by $|x \wedge y|$ and $|x \wedge y|^\diamond$ the Gromov product on $(X, {\mathcal E})$ and $(X, {\mathcal E}^\diamond)$ respectively.
 As a direct consequence of Theorem \ref{theorem-WSC-Hyperbolic} and \ref{theorem-quasi}, we have\\

\begin{theorem}
The graph $(X, {\mathcal E})$ is hyperbolic provided that the corresponding IFS satisfies the WSC or the self-similar set $K$ has  positive Lebesgue measure.
\end{theorem}
\begin{pf}
Observing that
\[
|x \wedge y| = \frac{1}{2}\big{(} |x| + |y| - d(x,y) \big{)}, \quad \mbox{and} \quad |x \wedge y|^\diamond = \frac{1}{2}\big{(} |x| + |y| - d^\diamond(x,y) \big{)}.
\]
It follows that
\begin{equation}\label{eq-deff-Gromov}
\Big{|} |x \wedge y| - |x \wedge y|^\diamond \Big{|} = \frac{1}{2} \Big{|} d(x,y) - d^\diamond(x,y) \Big{|} \le \frac{C}{2},
\end{equation}
where the constant $C>0$ is as in (\ref{eq-quasi}). \\

Note that $(X, {\mathcal E}^\diamond)$ is hyperbolic. Hence there exists a constant $\delta>0$ such that
\[
|x \wedge y|^\diamond \ge \min\{|x \wedge z|^\diamond, \ |z \wedge y|^\diamond\} - \delta, \quad \forall x,y,z \in X.
\]
Thus
\[
\begin{aligned}
|x \wedge y| & \ge |x \wedge y|^\diamond - \frac{C}{2}  \ge \min\{|x \wedge z|^\diamond, |z \wedge y|^\diamond\} - \delta - \frac{C}{2}\\
& \ge \min\{|x \wedge z|, |z \wedge y|\} - (\delta + C).
\end{aligned}
\]
This completes the proof.
\end{pf}

\bigskip

\begin{remark} \label{remark-section4}
In Theorem \ref{theorem-WSC-Hyperbolic} and \ref{theorem-quasi}, the IFS satisfies the WSC or the self-similar set has positive Lebesgue measure, and this implies the following condition:\\

\noindent {\rm
(C) For any $a>0$, there exists a constant $C>0$, such that for any integer $n > 0$, $D \subset {\mathbb R}^d, \ |D| \le a r^n$ and a subset $\{u_1, \cdots, u_\ell\} \subset {\mathcal J}_n$ with $S_{u_i}(K) \cap D \not=\emptyset$, $S_{u_i}(K) \cap S_{u_{i+1}}(K) \not= \emptyset$. Then there exist a subset $\{z_1, \cdots, z_{\ell_0}\} \subset \{u_1, \cdots, u_\ell\}$ such that $z_1 = u_1$, $z_{\ell_0} = u_\ell$, $S_{z_i}(K) \cap S_{z_{i+1}}(K) \not= \emptyset$ and $\ell_0 \le C$.\\
}
\indent
From the proof of Theorem \ref{theorem-WSC-Hyperbolic} and \ref{theorem-quasi}, we see that the above condition (C) implies both $(X, {\mathcal E})$ and $(X, {\mathcal E}^\diamond)$ are hyperbolic.
\end{remark}

\bigskip

\section{Hyperbolic Boundaries}

Throughout this section, we assume that the IFS satisfies the WSC or the self-similar set has positive Lebesgue measure, and hence the induced graphs $(X, {\mathcal E})$ and $(X, {\mathcal E}^\diamond)$ are hyperbolic. Denote by $\partial X$ and $\partial X^\diamond$ the hyperbolic boundaries, $\rho_a, \ \rho_a^\diamond$ the hyperbolic metrics with respectively.

It is know that if $f$  is a quasi-isometry from hyperbolic graph $(X_1, d_1)$ to $(X_2, d_2)$, then $\{x_n\}_n$ is Cauchy sequence in $X_1$ under the ultra-metric, if and only if $\{f(x_n)\}_n$ is. Moreover $\partial X_1$ and $\partial X_2$ are homeomorphism ( to see \cite{[CDP]} ). In our case, we have the following strengthen form. \\

\begin{prop}\label{theorem-equi-boundarys}
The hyperbolic boundaries $ \partial X = \partial X^\diamond$, and the hyperbolic metrics $\rho_a$ and $\rho_a^\diamond$ are equivalent, i.e., there exists a constant $C>0$ such that
\begin{equation}\label{eq-equi-boundary}
C^{-1} \rho_{a}( \xi, \eta) \le \rho_{a}^\diamond(\xi, \eta) \le C \rho_{a}(\xi, \eta), \quad \forall \xi, \eta \in \partial X.
\end{equation}
\end{prop}

\begin{pf}
Recall that a sequence $\{x_n\}_n \subset X$ with $|x_n| \to \infty$ is a Cauchy sequence under the ultra-metric $\rho_a$ if and only if $\lim_{m,n \to \infty} |x_m \wedge x_n| = \infty$, and the Cauchy sequence $\{y_n\}_n$ with $|x_n|=|y_n|$ equivalent to $\{x_n\}_n$ if and only if $d(x_n, y_n) \le c \delta$ for all but finite many $n$.

By (\ref{eq-deff-Gromov}), a sequence $\{x_n\}_n \subset X$ with $|x_n| \to \infty$ is a Cauchy sequence under the ultra-metric $\rho_{a}$ if and only if it is Cauchy in $\rho_{a}^\diamond$; moreover, by Theorem \ref{theorem-quasi}, the Cauchy sequence $\{y_n\}_n$ equivalent to $\{x_n\}_n$ in $\rho_{a}$ if and only if they are equivalent in $\rho_{a}^\diamond$. Recall an element in the hyperbolic boundary $\partial X$ is an equivalence class of Cauchy sequence in the ultra-metric $\rho_a$. Hence an element  $\xi \in \partial X$ if and only if $\xi \in \partial X^\diamond$, i.e., $\partial X = \partial X^\diamond$.

Now we prove (\ref{eq-equi-boundary}). For $\xi=[\{x_n\}_n], \ \eta=[\{y_n\}_n] \in \partial X = \partial X^\diamond$,  ($\xi \not= \eta$), by (\ref{eq-deff-Gromov}), we have
\[
\Big{|} |x_n \wedge y_n | - |x_n \wedge y_n|^\diamond \Big{|} \le \overline{C}, \quad \forall n.
\]
It follows that $\rho_{a}(x_n, y_n) \le e^{a \overline{C}} \rho_{a}^\diamond(x_n, y_n)$. Letting $n\to\infty$ and making use of (\ref{eq-ultra-ineq}), we have $\rho_{a}(\xi, \eta) \le C \rho_{a}^\diamond(\xi, \eta)$. The same argument implies the inverse inequality.
\end{pf}

\bigskip

To understand the topology of $(\partial X, \rho_{a})$ and $(\partial X^\diamond, \rho_{a}^\diamond)$, by the above Proposition, we need only to consider one of them. In the following, we consider $(\partial X, \rho_a)$. The arguments in \cite{[Lau-Wang-2009]} are adopted here.\\

\begin{lemma}\label{lemma-ray-II}
$\pi[u_0, u_1, \cdots]$ is a geodesic ray in the graph $(X, {\mathcal E})$ if and only if there exist $\bfi = i_1 i_2 \cdots \in \Sigma^\infty$ such that $u_n =[\bfi|_n]$ for all $n \ge 0$, where $\bfi|_n \in {\mathcal J}_n$ is the initial part of $\bfi$.
\end{lemma}
\begin{pf}
Clearly, for any $\bfi \in \Sigma^\infty$, $\pi[\bfi|_0, \bfi|_1, \cdots]$ is a geodesic ray in $(X, {\mathcal E})$ (where $\pi[\bfi|_0, \bfi|_1, \cdots]$ is abuse the notation for $\pi[[\bfi|_0], [\bfi|_1], \cdots]$). \\

Conversely, assume that $\pi[u_0, u_1, \cdots]$ is a geodesic ray in the graph $(X, {\mathcal E})$. Then for each $i\ge 0$, $(u_i, u_{i+1}) \in {\mathcal E}_0$. We use induction to construct $\bfi \in \Sigma^\infty$ as follow:

Choose any $\bfi_1 = i_1 i_2 \cdots i_k \in u_1$. If we have selected $\bfi_m = i_1, i_2 \cdots i_n \in u_\ell$, note that $(u_\ell, u_{\ell + 1}) \in {\mathcal E}_v$, by the definition of ${\mathcal E}_v$, we know that there are some $i_{n+1}, \cdots, i_{n+k'} \in \Sigma$ such that $i_1 \cdots i_n i_{n+1} \cdots i_{n+k'} \in u_{\ell + 1}$. Eventually we obtain the index $\bfi = i_1, i_2 \cdots \in \Sigma^\infty$ such that $\bfi|_n \in u_n$.
\end{pf}

\begin{lemma} \label{lemma-equi-ray}
Let $\xi = \pi[u_0, u_1, \cdots]$ be a geodesic ray in $(X, {\mathcal E})$. Then the limit $\lim_{n\to \infty} S_{u_n}(x)$ exists and is independent of $x\in {\mathbb R}^d$. Moreover, if a geodesic ray $\eta=\pi[v_0, v_1, \cdots]$ is equivalent to $\xi$, then $\lim_{n\to\infty} S_{v_n}(x) = \lim_{n\to\infty} S_{u_n}(x)$.
\end{lemma}
\begin{pf}
By Lemma \ref{lemma-ray-II}, there exists $\bfi \in \Sigma^\infty$ such that $u_n = [\bfi|_n]$, $(n=0,1, \cdots)$. It is well known that the limit $\lim_{n\to\infty} S_{\bfi|_n}(x) = \lim_{n\to\infty} S_{u_n}(x)$ exists and is independent of $x \in {\mathbb R}^d$, and the first part of the lemma follows.

For the second part, note that $\xi$ and $\eta$ are equivalent, and hence there exists some constant $c>0$ such that
$d(u_n, v_n) \le c$ for all $n \ge 0$. For each fixed $n$, there is a geodesic segment $\pi[w_1, \cdots, w_\ell]$ ($w_1 = u_n, \ w_\ell = v_n$ and $\ell \le c$) connecting  $u_n$ and $v_n$. Note that $S_{w_i}(K) \cap S_{w_{i+1}}(K) \not= \emptyset$, taking any $x\in K$, we have the following estimate
\[
| S_{u_n}(x) - S_{v_n}(x)| \le \sum_{i=1}^{\ell} {\rm diam} S_{w_i}(K) \le \ell |K| r^{n- \ell} \le c |K| r^{n-c}.
\]
This implies that $\lim_{n\to\infty} S_{v_n}(x) = \lim_{n\to\infty} S_{u_n}(x)$.
\end{pf}

\bigskip

Let $\xi = \pi[u_0, u_1, \cdots]$ be a geodesic ray in $(X, {\mathcal E})$. We define
\[
\Phi(\xi) = \lim_{n\to\infty}S_{u_n}(x_0),
\]
where $x_0\in {\mathbb R}^d$. By using the above lemma, if a geodesic ray $\eta$ is equivalent to $\xi$, then $\Phi(\xi) = \Phi(\eta)$. Hence $\Phi$ induces a map (we still use $\Phi$ to denote this map) from the hyperbolic boundary $\partial X$ to the self-similar set $K$.\\

\begin{theorem} \label{theorem-homeomorphism}
The map $\Phi:\ \partial X \longrightarrow K$ is a bijection and there exists a constant $C>0$ such that
\begin{equation}\label{eq-holder-1}
|\Phi(\xi) - \Phi(\eta)| \le C \rho_{a}(\xi,\eta)^\alpha, \quad \forall \xi, \eta \in \partial X,
\end{equation}
where  $\alpha= - \log r/a$. In particular $\partial X$ is homeomorphic to the self-similar set $K$.
\end{theorem}

\begin{pf}
Let $x_0 \in {\mathbb R}^d$. For any $x\in K$, there exists an index $\bfu=i_1 i_2 \cdots \in \Sigma^\infty$ such that $\lim_{n\to\infty} S_{i_1 i_2 \cdots i_n} (x_0) = x$. In particular, $\lim_{n\to\infty}S_{\bfu|_n}(x_0) = x$ (recall that $\bfu|_n = i_1 \cdots i_k$, where the integer $k$ is such that $r_1 r_2 \cdots r_k \le r^n < r_1 r_2 \cdots r^{k-1}$). This means that the image of the geodesic ray $\pi[\bfu|_0, \bfu|_1, \bfu|_2, \cdots]$ under the map $\Phi$ is $x\in K$. Hence the map is surjective. \\

To show that $\Phi$ is injective, assume that $\xi, \eta \in \partial X$. Then there are geodesic rays $\pi[x_0, x_1, \cdots]$ and $\pi[y_0, y_1, \cdots]$ converge to $\xi$ and $\eta$  respectively; moreover we assume that they attain  the infimum  in (\ref{eq-extend-Gromov}).

By Lemma \ref{lemma-ray-II}, there exist indexes $\bfu =i_1 i_2 \cdots , \ \bfv =j_1 j_2 \cdots \in \Sigma^\infty$ such that  $x_n = \bfu|_n$ and $y_n = \bfv|_n, \ n=0, 1, \cdots$ (recall that we abuse the notation $x_n = \bfu|_n$ means that $ x_n = [\bfu|_n]$). If  $\Phi(\xi) = \Phi(\eta) = x \in K$, then
\[
x \in S_{x_n}(K) \bigcap S_{y_n}(K), \ \quad  n =0 ,1, 2, \cdots.
\]
Hence $(x_n, y_n) \in {\mathcal E}$. It follows that $d(x_n, y_n) \le 1$, we see that the geodesic rays $\pi[x_0, x_1, \cdots]$ and $\pi[y_0, y_1, \cdots]$ are equivalent, i.e., $\xi = \eta$. Hence the map $\Phi$ is injective. \\

Now we prove (\ref{eq-holder-1}). If the above $\xi \not= \eta$, then for any fixed $n$, there is a canonical geodesic $\pi[z_{0,n}, z_{1,n}, \cdots, z_{k_n, n}]$ ($z_{0,n}= x_n, \ z_{k_n, n} = y_n$) joining $x_n$ and $y_n$. Note that $x_{n+1}, z_{0,n} \cdots, z_{k_n, n}, y_{n+1}$ is a path (may not geodesic) from $x_{n+1}$ to $y_{n+1}$. Hence $d(x_{n+1}, y_{n+1}) \le d(x_n, y_n) + 2$. It follows that
\[
\begin{aligned}
| x_{n+1} \wedge y_{n+1} | & = \frac{1}{2} ( |x_{n+1}| + |y_{n+1}| - d(x_{n+1}, y_{n+1}) )\\
& \ge \frac{1}{2} ( |x_n| + |y_n| - d(x_{n}, y_{n}) ) = | x_n \wedge y_n |.
\end{aligned}
\]
i.e., $\{| x_n \wedge y_n | \}_{n=1}^\infty$ is an increasing  sequence. On the other hand, $\xi \not= \eta$  implies that
\[
 |\xi \wedge \eta | = \lim_{n\to\infty} | x_n \wedge y_n | = k \ < \infty.
\]
Note that $2 k$ is an integer. Hence there exists $m$ such that
\[
| x_n \wedge y_n |  = k, \mbox{ if } n \ge m; \quad  | x_n \wedge y_n | < k,  \mbox{ if } n < m.
\]

To estimate $|\Phi(\xi) - \Phi(\eta)|$, we note that $\Phi(\xi) \in S_{x_{m}}(K)$ and $\Phi(\eta) \in S_{y_{m}}(K)$, hence there exist $x,y \in K$ such that $S_{x_{m}}(x)= \Phi(\xi)$ and $S_{y_{m}}(y) = \Phi(\eta)$. Recall that $\pi[z_{0,m}, \cdots, z_{k_m, m}]$ is a canonical geodesic jointing $x_m$ and $y_m$.
Assume that $z_{i, m}, \cdots, z_{i', m}$, $(0 \le i \le i' \le k_m)$ is the horizontal part of the canonical geodesic, denote by $k', \ \ell(=i'-i)$ the level and length of this segment. Then by (\ref{eq-graph2-Gromov}), we have
\[
| x_m \wedge y_m | = k = k' + \frac{\ell}{2}.
\]
For the vertical parts, we have
\[
\Phi(\xi) = S_{x_m}(x) \in  S_{x_m}(K) \subset S_{z_{i,m}}(K), \quad
\Phi(\eta) = S_{y_m}(x) \in S_{y_m}(K) \subset S_{z_{i',m}}(K).
\]
It follows that
\[
| \Phi(\xi) - \Phi(\eta)|  \le {\rm diam} \Big{(} \bigcup_{j=i}^{i'} S_{z_{j,m}}(K) \Big{)}
\le \sum_{j=i}^{i'} {\rm diam} \Big{(} S_{z_{j,m}}(K) \Big{)}
\le (\ell+1) |K| r^{k'}
\]
By Theorem \ref{theorem-graph-II-hyperbloic}, $\ell$ is bounded by the constant $L>0$. Hence
\[
| \Phi(\xi) - \Phi(\eta)|  \le (L+1) |K| r^{k- \frac{\ell}{2}} \le (L+1) r^{-\frac{L}{2}} |K| \rho_{a,2}(\xi, \eta)^\alpha,
\]
where $\alpha = -\log r/a$, and (\ref{eq-holder-1}) follows.

By (\ref{eq-holder-1}), we know that the map $\Phi$ is continuous, and hence is a homeomorphism, since $\partial X$ and $K$ are compact. We complete the proof.
\end{pf}

\bigskip

In order to get the inverse inequality of (\ref{eq-holder-1}), we need the following condition on IFS as \cite{[Lau-Wang-2009]}:\\

\noindent
{\it
(H) There exists a constant $C'>0$ such that for any integer $n>0$ and $\bfu, \bfv \in {\mathcal J}_n$, either \[
S_\bfu(K) \cap S_\bfv(K) \not= \emptyset \quad \mbox{ or } \quad |S_\bfu(x) - S_\bfv(y)| \ge C' r^n, \quad \forall x,y\in K.
\]
}

\bigskip

\begin{prop}\label{prop-holder}
Suppose the IFS $\{S_j\}_{j=1}^N$ in Theorem \ref{theorem-homeomorphism} satisfies in addition condition (H). Then there exists a constant $C>0$ such that for any $\xi, \eta \in \partial X$,
\begin{equation}\label{eq-holder-2}
C^{-1} | \Phi(\xi) - \Phi(\eta) | \le \rho_{a}(\xi, \eta)^\alpha \le C | \Phi(\xi) - \Phi(\eta) |,
\end{equation}
where $\alpha= -\log r /a$.
\end{prop}
\begin{pf}
For $\xi=\eta$ the inequality is trivial, in the following we assume that $\xi \not= \eta$ and use the notation  in the proof of Theorem \ref{theorem-homeomorphism}. It is clear that $S_{x_{k' + 1}}(K) \cap S_{y_{k'+1}}(K) = \emptyset$, where $k'$ is the level of the horizontal part of the canonical geodesic connecting $\xi$ and $\eta$ as in proof of the above Theorem.

By condition (H), we have
\[
|\Phi(\xi) - \Phi(\eta)| \ge C' r^{k' + 1} = C' r^{1- \frac{\ell}{2}} r^k \ge C \rho_a(\xi,\eta)^\alpha
\]
for some constant $C > 0$. This is the second inequality of (\ref{eq-holder-2}). The first inequality is proved in Theorem \ref{theorem-homeomorphism}.
\end{pf}

\bigskip

The above theorem can be used to study the Lipschitz equivalence relation for self-similar sets which can be found in \cite{[LL]}. \\

The following example shows that the second inequality of (\ref{eq-holder-2}) does not hold if the condition (H) fails. \\

\noindent
{\bf Example 2:} Let $p_0=(0,0),\ p_1=(1,0)$ and $p_2=(\frac{1}{2}, \frac{\sqrt{3}}{2})$ be the vertices of equilateral triangle $\Delta$ in ${\mathbb R}^2$, and let $S_i(x) = \frac{1}{3}(x+q_i), \ q_i \in {\mathbb R}^2, \ i= 0,1, \cdots, 4$ be the five maps, each maping the triangle $\Delta$ to a small triangle (to see Fig 3}). This IFS satisfies the OSC with the interior of $\Delta$ as an open set. Denote by $K$ the self-similar set of the IFS $\{S_i\}_{i=0}^4$. Then $\{p_0, p_1, p_2\} \subset K \subset \Delta$.
Let $x_0$ be the horizontal coordinate of left-bottom of triangle $S_3(\Delta)$ (the first coordinate of $S_3(p_0)$), we choose $x_0$ such that
\[
\frac{1}{2} - x_0 = \sum_{k=1}^\infty 3^{-n_k}, \quad n_k = 1 + \frac{k(k+1)}{2}, \ k =1,2, \cdots.
\]
In this example $r=\frac{1}{3}, \ X = \Sigma^* = \cup_{n=0}^\infty \{0,1,2,3,4\}^n$.

We define a sequence $\{a_i\}_{i=1}^\infty$  in the symbolic space $\Sigma = \{0,1, \cdots, 4\}$ as: $a_i=1$ if $i=n_k$ for some $k$; otherwise $a_i=0$. By using this sequence, we prove that condition (H) does not hold.\\

\begin{center}
\begin{figure}[h]  \label{fig4}
\centerline{\includegraphics[width=5cm,height=4cm]{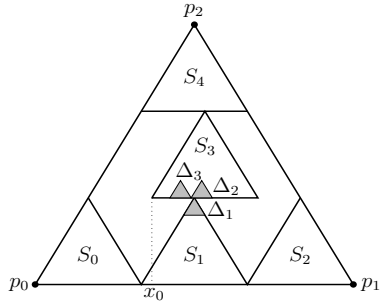}}
\caption{\small{The maps in Example 2. }}
\end{figure}
\end{center}

We iterate the IFS $n_k$ times and get a set of small triangles $\{S_\bfu (\Delta):\ |\bfu|= n_k\}$. Consider the following three of such trangles: $\Delta_1$ is the one on the top of $S_1(\Delta)$, $\Delta_2$ the unique one in $S_3(\Delta)$ which intersects $\Delta_1$, and $\Delta_3$ the one on the left of $\Delta_2$. The corresponding codes are  $\bfu_k=1 4^{n_k-1}$, $\bfv_k = 3 a_2 a_3 \cdots a_{n_k}$ and $\bfw_k = 3 a_2 a_3 \cdots a_{n_k-1}0$ with respectively. Let $x_k$ be the  coordinate of right-bottom of triangle $\Delta_3$, i.e., $x_k= S_{\bfw_k}(p_1)$, and let $y= S_{\bfu_k}(p_2)$. Then
\[
| y- x_k | = \min\{ |x-x'|:\ x \in \Delta_1,\ x'\in \Delta_3\} = \sum_{i=k+1}^\infty 3^{-n_i} = c_k \cdot 3^{-n_{k+1}},
\]
where $1< c_k < \frac{3}{2}$. We see that the condition (H) does not hold.\\

Consider the geodesic rays $\xi_k= \{\bfw_k 2^\infty|_n\}_n$ and $\eta= \{1 4^\infty|_n\}$. Then $\Phi(\xi_k) = x_k$ and $\Phi(\eta) = y$. On the other hand,
\[
|\xi_k \wedge \eta| = |\bfw_k| + 1 =n_k + 1,
\]
and $\rho_{a}(\xi_k, \eta)^\alpha = 3^{- n_k - 1}$. We see that the second inequality in (\ref{eq-holder-2}) does not hold. {\hfill$\Box$}

\bigskip
Another example which does not satisfy the condition (H) can be found in \cite{[Lau-Wang-2009]}.

\bigskip

To end this paper, we ask the following question:\\
{\bf Question 1:} Does that graph $(X, {\mathcal E})$ or $(X, {\mathcal E}^\diamond)$ is local finite imply that the IFS satisfies the WSC?\\
{\bf Question 2:} Are the conditions that the IFS satisfies the WSC or the self-similar set has  positive Lebesgue measure necessary for the graphs to be hyperbolic?\\

\vspace{2ex} \noindent {\bf Acknowledgements.}~~  Part of this work
was carried out while the  author was visiting the
Department of Mathematics of the Chinese University of Hong Kong.
He is grateful for the discussion with professors K.S. Lau and S.M. Ngai.

\end{document}